\documentclass[12pt,a4paper,draft]{article}
\usepackage[T2A]{fontenc}
\usepackage[english]{babel}
\usepackage{latexsym,amsfonts,amssymb,amsmath,longtable,amsthm,mathrsfs}
\usepackage[unicode,final,hyperindex]{hyperref}
\renewcommand{\le}{\leqslant}
\renewcommand{\ge}{\geqslant}

\newtheorem{Lemma}{{\bfseries Lemma}}
\newtheorem{Cor}[Lemma]{{\bfseries Corollary}}
\newtheorem{Theo}[Lemma]{{\bfseries Theorem}}
\theoremstyle{definition}

\title{\vspace{-1cm} \hfill{\normalsize 20D20}{
\fontfamily{cmr} \fontseries{bx} \selectfont \\ \vspace{1cm} On the intersection of solvable Hall subgroups \\ in finite simple exceptional \\ groups of Lie type}\footnote{The work is supported by RFBR, projects 11--01--00456, 12--01--33102.}}
\date{}
\author{\bf   Evgeny. P. Vdovin}

\begin{document}

\maketitle

\begin{abstract}
Assume that a finite almost simple group with simple socle isomorphic to an exceptional group of Lie type possesses a solvable Hall subgroup. Then there exist four conjugates of the subgroup such that their intersection is trivial.


\end{abstract}

\section*{Introduction}

Throughout the paper the term ``group'' we always use in the meaning ``finite group''.  We use symbols $A\leq G$ and
$A\unlhd G$ if $A$ is a subgroup of $G$ and  $A$ is a normal subgroup of $G$ respectively.
If $\Omega$ is a (finite) set, then by $\mathrm{Sym}(\Omega)$ we denote the group of all permutations of
$\Omega$. We also denote $\mathrm{Sym}(\{1,\ldots,n\})$ by $\mathrm{Sym}_n$. Given $H\leq G$ by
$H_G=\cap_{g\in G} H^g$ we denote the kernel of~$H$.

Assume that $G$ acts on $\Omega$. An element $x\in
\Omega$ is called a {\em $G$-regular point}, if $\vert
xG\vert=\vert G\vert$, i.e., if the stabilizer of $x$ is trivial. We define the action of $G$ on $\Omega^k$
by
\begin{equation*}
g:(i_1,\ldots,i_k)\mapsto
(i_1g,\ldots,i_kg).
\end{equation*}
If $G$ acts faithfully and transitively on $\Omega$, then the minimal $k$ such that
$\Omega^k$ possesses a $G$-regular point is called the {\em base size} of $G$ and is denoted by~$\mathrm{Base}(G)$.
For every natural $m$ the number of $G$-regular orbits on  $\Omega^m$ is denoted by
$\mathrm{Reg}(G,m)$ (this number equals $0$ if $m<\mathrm{Base}(G)$). If $H$ is a subgroup of $G$ and $G$
acts on the set $\Omega$ of right cosets of $H$ by right multiplications, then
$G/H_G$ acts faithfully and transitively on $\Omega$. In this case we denote  $\mathrm{Base}(G/H_G)$
and $\mathrm{Reg}(G/H_G,m)$ by $\mathrm{Base}_H(G)$ and $\mathrm{Reg}_H(G,m)$ respectively. We also say that
$\mathrm{Base}_H(G)$ is the {\em base size of}  $G$ {\em with respect to}~$H$. Clearly,
$\mathrm{Base}_H(G)$ is the minimal $k$ such that there exist elements $x_1,\ldots,x_k\in G$ with
$H^{x_1}\cap \ldots\cap H^{x_k}=H_G$. Thus, the base size of $G$ with respect to $H$ is the minimal $k$ such that there exist $k$ conjugates of $H$ with intersection equals~$H_G$.

We prove the following theorem in the  paper.

\begin{Theo}\label{IntersectSimple}
{\em (Main Theorem)} Let $G$ be an almost simple group with simple socle isomorphic to an exceptional group of Lie type. Assume also that $G$ possesses a solvable Hall subgroup $H$. Then  $\mathrm{Base}_H(G)\le 4$.
\end{Theo}

The following results were obtained in this direction. In 1966 D.S.Passman proved (see \cite{Pass}) that a
$p$-solvable group possesses three Sylow $p$-subgroups whose intersection equals the $p$-radical of $G$. Later in 1996 V.I.Zenkov proved (see \cite{Zen1}) that the same conclusion holds for arbitrary finite group $G$. In \cite{Dolfi} S.Dolfi proved that in every $\pi$-solvable group $G$
there exist three conjugate $\pi$-Hall subgroups whose intersection equals $O_\pi(G)$ (see also~\cite{VdovinIntersSolv}). Notice also that
V.I.Zenkov in \cite{Zen2} constructed an example of a group $G$ possessing a solvable $\pi$-Hall subgroup
$H$ such that the intersection of five conjugates of $H$ equals $O_\pi(G)$, while the intersection of every four conjugates of  $H$
is greater than~$O_\pi(G)$. In
\cite[Theorem~1]{VdovinZenkov} the following statement is proven.

\begin{Theo}\label{Vdo-Zen-main}
Let $G$ be a finite group possessing a solvable $\pi$-Hall subgroup $H$. Assume that for every simple component $S$ of $E(\overline{G})$
of the factor group $\overline{G}=G/S(G)$, where $S(G)$ is a solvable radical of $G$, the following condition holds:
\begin{multline*}
\text{for every  } L\text{ such that } S\leq L\leq \mathrm{Aut(S)}\text{ and contains a solvable }\pi-\text{Hall subgroup
}M,\\ \text{ the inequalities }\mathrm{Base}_M(L)\le 5\text{ and }\mathrm{Reg}_M(L,5)\ge5\text{ hold.}
\end{multline*}
Then $\mathrm{Base}_H(G)\le 5$ and $\mathrm{Reg}_H(G,5)\ge5$.
\end{Theo}

Moreover, at the beginning of the proof of Theorem 2 from \cite{VdovinZenkov} the following statement is obtained.

\begin{Lemma}\label{BoundReg}
If, for a group $G$ and its subgroup $H$, the inequality $\mathrm{Base}_H(G)\leq4$ holds, then
$\mathrm{Reg}_H(G,5)\ge5$.
\end{Lemma}

Thus by Theorems \ref{IntersectSimple} and \ref{Vdo-Zen-main}, Lemma \ref{BoundReg}, and  \cite[Theorem~2]{VdovinZenkov} we immediately obtain

\begin{Theo}\label{IntersectFinite}
Let $H$ be a solvable $\pi$-Hall subgroup of $G$. Assume that each nonabelian composition factor of the socle of
$G/S(G)$, where $S(G)$ is the solvable radical of $G$, is isomorphic to either alternative, or sporadic, or exceptional group of Lie type. Then
$\mathrm{Base}_H(G)\le 5$, i.e., there exist elements $x,y,z,t$ of $G$ such that the identity $$H\cap H^x\cap H^y\cap H^z\cap H^t=O_\pi(G)$$ holds.
\end{Theo}

\section{Notations and preliminary results}       

Throughout by $\pi$ a set of primes is denoted, while by $\pi'$ we denote its complement in the set of all primes.
A subgroup $H$ of $G$ is called a {\em $\pi$-Hall subgroup}, if the order $\vert H\vert$ is divisible by primes from $\pi$ only,
while its index $\vert G:H\vert$ is divisible by primes from $\pi'$ only. The set of all $\pi$-Hall subgroups of $G$ is denoted by $\mathrm{Hall}_\pi(G)$. A subgroup $H$ of $G$ is called a {\em Hall subgroup}, if its order $\vert H\vert$ and the index $\vert G:H\vert$ are coprime. A group $G$ is called {\em almost simple}, if there exists a nonabelian simple group $S$ such that $F^\ast(G)=S$, where $F^\ast(G)$ is the generalized Fitting subgroup of~$G$. In other words, $G$ is called almost simple, if there exists a simple group $S$ such that $S\simeq \mathrm{Inn}(S)\leq G\leq \mathrm{Aut}(S)$.

\begin{Lemma}\label{HallSubgrprop} {\em \cite[Lemma~1]{Hall}} Let  $G$ be a finite group and $A$ be its normal subgroup.
If $H\in \mathrm{ Hall}_\pi(G)$, then  ${H\cap A\in \mathrm{ Hall}_\pi(A)}$ and
${HA/A\in \mathrm{Hall}_\pi(G/A)}$.
\end{Lemma}

\begin{Lemma}\label{ZenkovAbelian} {\em \cite{Zen3}}
Let $A$ be an abelian subgroup of a finite group $G$. Then there exists $x\in G$ such that $A\cap A^x \leq F(G)$.
\end{Lemma}

Combining known results (see \cite[Theorems 8.3--8.7]{VdoRevUspehi}), we obtain the following

\begin{Lemma}\label{SolvHallcharacteristic}
Let $G$ be a simple group of Lie type over a field of characteristic $p\in \pi$ and $H$ be its solvable $\pi$-Hall subgroup. Then either
$H$ is included in a Borel subgroup of $G$, or one of the following holds:
\begin{itemize}
\item[{\em (1)}] $G=SL_3(2)$ or $G=SL_3(3)$ and $H$ is the stabilizer of a line or of a plain in the natural
$3$-dimensional module, i.e., there exist two classes of conjugate $\pi$-Hall subgroups in this case.
\item[{\em (2)}] $G=SL_4(2)$ or $G=PSL_4(3)$ and $H$ is the stabilizer of a two-dimensional subspace of the natural $4$-dimensional module.
\item[{\em (3)}] $G=SL_5(2)$ or $G=SL_5(3)$ and $H$ is the stabilizer of a chain of subspaces
$V_0<V_1<V_2<V_3=V$ whose codimensions are in the set $\{1,2\}$ (i.e. two codimensions equal $2$ and one codimension equals $1$). There exist three classes of conjugate $\pi$-Hall subgroups in this case.
\end{itemize}
\end{Lemma}

We recall some known technical results (see \cite{BurLieSha}). If $G$ acts transitively on the set $\Omega$, then given $x\in G$ by $\mathrm{fpr}(x)$ we denote the fixed point ratio of $x$, i.e. $\mathrm{fpr}(x)=\vert \mathrm{fix}(x)\vert/\vert
\Omega\vert$, where $\mathrm{fix}(x)=\{\omega\in\Omega\mid \omega^x=\omega\}$. If $G$ acts transitively and
$H$ is a point stabilizer, then the following formulae is known
\begin{equation}\label{fpr}
\mathrm{fpr}(x)=\frac{\vert x^G\cap H\vert}{\vert x^G\vert}.
\end{equation}
As it is noted in \cite[Theorem~1.3]{LieSha}, the base size can be bounded by using the following arguments. Assume that $G$ acts faithfully and let $Q(G,c)$ denote the probability that arbitrary chosen element of $\Omega^c$ is not a $G$-regular point. Clearly,
$\mathrm{Base}(G)$ is the minimal $c$ such that $Q(G,c)<1$. In particular, if $Q(G,c)<1$ then $\mathrm{Base}(G)\le
c$. Clearly, an element of $\Omega^c$ is not a $G$-regular point if and only if it is stable under the action of an element
$x$ of prime order. Notice also that the probability for arbitrary chosen element of $\Omega^c$ to be stable under $x$ is not greater than
$\mathrm{fpr}(x)^c$. Denote by $\mathcal{P}$ the set of elements of $G$ whose order is  equal to a prime number. Let $x_1,\ldots,x_k$ be representatives of the conjugacy classes of elements from $\mathcal{P}$. Since
$G$ acts transitively, the formulae  \eqref{fpr} shows that $\mathrm{fpr}(x)$ does not depend on the choice of the representative of a conjugacy class. Thus the following chain of inequalities holds.
\begin{equation}\label{BoundBase}
Q(G,c)\le \sum_{x\in\mathcal{P}} \mathrm{fpr}(x)^c=\sum_{i=1}^k\vert x_i^G\vert\cdot
\mathrm{fpr}(x_i)^c=:\widehat{Q}(G,c).
\end{equation}

In particular, we can use the upper bound for $\mathrm{fpr}(x)$ in order to bound $\widehat{Q}(G,c)$ and so to bound $Q(G,c)$. The following lemma is the main technical tool for this bound.

\begin{Lemma}\label{BoundfprBase}
{\em \cite[Proposition~2.3]{BurLieSha}} Let $G$ be a transitive group of permutations on $\Omega$ and $H$ be a point stabilizer.
Assume that $x_1,\ldots,x_k$ are representatives of distinct conjugacy classes such that the inequalities $\sum_i\vert x_i^G\cap H\vert\le A$ and $\vert x_i^G\vert\ge B$ hold for all $i=1,\ldots,k$.  Then the inequality
\begin{equation*}
\sum_{i=1}^k\vert x_i^G\vert\cdot
\mathrm{fpr}(x_i)^c\leq B(A/B)^c
\end{equation*}
holds for every $c\in \mathbb{N}$.
\end{Lemma}

Notice that for every subgroup $H$ and every set $x_1,\ldots,x_k$ not containing the identity element the bound $\sum_i\vert x_i^G\cap H\vert< \vert H\vert$ holds.

\section{Technical results}

Our notations for groups of Lie type agree with that of \cite{GorLySol}. In particular, for every simple group of Lie type $S$
over a field of characteristic $p$ we fix a simple algebraic group $\overline{G}$ of adjoint type and a Steinberg map $\sigma$
so that $S=O^{p'}(\overline{G}_\sigma)$. Then $\overline{G}_\sigma$ is the group of inner-diagonal  automorphisms of
$S$ (we denote the group of inner-diagonal automorphisms of $S$ by $\widehat{S}$). We assume that a Borel $\overline{B}$ and its maximal torus $\overline{T}$ are chosen $\sigma$-invariant, and we denote $\overline{B}_\sigma$ and $\overline{T}_\sigma$ by  $B$ and $T$ respectively. Recall
that if  $S\in \{{}^2A_n(q),{}^2D_n(q),{}^2E_6(q)\}$, then the definition field of $S$ equals $\mathbb{F}_{q^2}$, if
$S={}^3D_4(q)$, then the definition field of $S$ equals  $\mathbb{F}_{q^3}$, and the definition field of $S$ equals
$\mathbb{F}_{q}$ in the remaining cases. For groups  ${}^2A_n(q),{}^2D_n(q),{}^2E_6(q)$ we also use the notations $A_n^-(q),D_n^-(q),E_6^-(q)$ respectively. Notice also the known fact: $Z(\overline{B})\cap
\overline{T}=Z(\overline{G})$ ($=1$, if $\overline{G}$ is of adjoint type) and $Z(B)\cap T=Z(S)$ ($=1$, if $\overline{G}$ is of adjoint type).

\begin{Lemma}\label{IntersectionParabolicLeviCenter}
Let $G$ be a group of inner-diagonal automorphisms of a finite simple group of Lie type over a field of characteristic
$p$ (i.e. $G=\overline{G}_\sigma$ for some connected simple algebraic group $\overline{G}$ of adjoint type
over an algebraically closed field of characteristic $p$ and a Steinberg map $\sigma$). Let
$B=U\leftthreetimes T$ be a Borel subgroup of $G$, where $U$ is a maximal unipotent subgroup of  $G$ and $T$ is a Cartan subgroup of $G$.
We denote the subgroup of monomial matrices containing  $T$ by $N$ so that $N/T\simeq W$ is the Weyl group of $G$.
Let $w_0\in W$ be the unique element that maps all positive roots into negatives, and  $n_0$ be its preimage in~$N$. Then there exists $x\in U^{n_0}$ such that $T^x\cap B=1$. In particular, there exist
$u,v\in O^{p'}(G)$ such that $B\cap B^u\cap B^v=1$.
\end{Lemma}

\begin{proof}
Consider $B^{n_0}=U^{n_0}\leftthreetimes T$. The Fitting subgroup $F(U^{n_0}\leftthreetimes
T)$ equals $U^{n_0}$ since $Z(O^{p'}(G))=1$. Otherwise, since $U^{n_0}$ is a normal nilpotent subgroup of
$U^{n_0}\leftthreetimes T$ we obtain that $U^{n_0}\leq F(U^{n_0}\leftthreetimes T)$. If $U^{n_0}\neq
F(U^{n_0}\leftthreetimes T)$, then there exists $1\neq z\in T$ centralizing $U^{n_0}$ and so lying in
$Z(O^{p'}(G))=1$, a contradiction. Hence $F(U^{n_0}\leftthreetimes
T)=U^{n_0}$ and by Lemma \ref{ZenkovAbelian} there exists $x\in U^{n_0}$ such that
$T\cap T^x=1$.

Notice that  $U^{n_0}\cap B=1$, so $(U^{n_0}\leftthreetimes T)\cap B=T$. Since $T^x\in U^{n_0}\leftthreetimes T$
we obtain
$$1=T^x\cap T=T^x\cap \left(\left(U^{n_0}\leftthreetimes T\right)\cap B\right)=\left(T^x\cap
\left(U^{n_0}\leftthreetimes T\right)\right)\cap B=T^x\cap B,$$ whence the main statement of the lemma follows.

Now we prove ``in particular'', i.e., we show that there exist $u,v\in O^{p'}(G)$ such that $B\cap B^u\cap
B^v=1$. By construction,  $x\in U^{n_0}\leq O^{p'}(G)$ and $1=T^x\cap B=(B^{n_0}\cap
B)^x\cap B=B\cap B^x\cap B^{n_0x}$. The lemma is proven.
\end{proof}

Let $S=O^{p'}(\overline{G}_\sigma)$ be a finite simple nontwisted group of Lie type over a field $\mathbb{F}_q$ of characteristic $p$.
A Cartan subgroup $T\cap S$ of $S$ can be obtained as $\langle h_r(\lambda)\mid r\in\Pi,
\lambda\in\mathbb{F}_q^\ast\rangle$ (see \cite[Theorem~2.4.7]{GorLySol}), where $\Pi$ is a set of fundamental roots of the root system of $S$.
Then a field automorphism $\varphi$ of $S$ can be chosen so that for every $r\in\Pi$, $\lambda\in\mathbb{F}_q^\ast$ the identity
$h_r(\lambda)^\varphi=h_r(\lambda^p)$ holds. Moreover, a graph automorphism  $\tau$ corresponding to the symmetry  $\overline{\phantom{\Pi}}:\Pi\rightarrow \Pi$ of the Dynkin diagram of $S$ can be chosen so that for every
$r\in\Pi$, $\lambda\in\mathbb{F}_q^\ast$ the identity $(h_r(\lambda))^\tau =h_{\bar{r}}(\bar{\lambda})$ holds, where
$\bar{\lambda}=\lambda$, if all roots have the same length. Consider the subgroup $A$ generated by so chosen field automorphism and graph automorphisms  (there exist several graph automorphisms for the root system $D_4$). It is well-known that $\mathrm{Aut}(S)=\widehat{S}\leftthreetimes A$. Moreover, $A$ normalizes a Borel subgroup $B$ containing the Cartan subgroup $T$. Since $N_{\widehat{S}}(B)=B$ we obtain that~${N_{\mathrm{Aut}(S)}(B)=B\leftthreetimes A}$.

Now assume that $S$ is a finite simple twisted group of Lie type distinct from a Suzuki group or a Ree group,  $L$ is a nontwisted group of Lie type
and $\psi$ is an automorphism of $L$ such that  $S=O^{p'}(L_\psi)$. Let
$\overline{\phantom{\Pi}}:\Pi\rightarrow \Pi$ be the symmetry of the Dynkin diagram of a fundamental set of roots $\Pi$ of the root system of $L$ using for construction of~$\psi$. Then a Cartan subgroup $T\cap L$ of $\widehat{L}$ can be written as $\langle h_r(\lambda)\mid r\in\Pi, \lambda\in\mathbb{F}_q^\ast\rangle$, and a field automorphism
$\varphi$ of $S$ can be chosen so that  for every $r\in\Pi$, $\lambda\in\mathbb{F}_q^\ast$ the equality
$(h_r(\lambda))^\varphi =h_{\bar{r}}(\lambda^p)$ holds. We set $A=\langle\varphi\rangle$, then
$\mathrm{Aut}(S)=\widehat{S}\leftthreetimes A$, and there exists a Borel subgroup $B$ of $\widehat{S}$ such that the equality~${N_{\mathrm{Aut}(S)}(B)=B\leftthreetimes A}$ holds.

\begin{Lemma}\label{IntersectionOutBorel}
In the introduced notations assume that, if  $S$ is not twisted, then the order $q$ of the definition field $\mathbb{F}_q$ of $S$ is greater than $2$. Moreover, if $S=D_4(q)$, assume also that $q>3$. Assume also that $S$ is neither a Suzuki group nor a Ree group. Then there exists $x\in T\cap S$ such that $C_A(x)=1$. In particular $A\cap A^x=1$.
\end{Lemma}

\begin{proof} If $S$ is not twisted and is distinct from  $D_4(q)$, then we can take $x=h_r(\lambda)$,
where $r\in\Pi$ is such that $r\not=\bar{r}$ and $\lambda$ is a generating element of the multiplicative group of
$\mathbb{F}_q$. If  $S$ is twisted distinct from  ${}^3D_4(q)$, then we can take $x=h_r(\lambda)h_{\bar{r}}(\lambda^q)$, where $\lambda$ is a   generating element of the multiplicative group of
$\mathbb{F}_{q^2}$ and $r\not=\bar{r}$. If $S={}^3D_4(q)$, then we can take
$x=h_r(\lambda)h_{\bar{r}}(\lambda^q)h_{\bar{\bar{r}}}(\lambda^{q^2})$, where $\lambda$ is a generating element of the multiplicative group of
$\mathbb{F}_{q^3}$ and $r\not=\bar{r}$. Finally, if $S=D_4(q)$ and $q>3$, then there exist $\lambda_1,\lambda_2\in\mathbb{F}_q^\ast\setminus\{1\}$ such that $\lambda_2\not\in\{\lambda_1^p,\lambda_1^{p^2},\ldots,\lambda_1^{q}\}$ and $\lambda_1$ generates $\mathbb{F}_q^\ast$. Choose fundamental roots $r,s$ so that there exists a nontrivial symmetry of the Dynkin diagram, permuting the roots. Then we can take $x=h_r(\lambda_1)h_s(\lambda_2)$.
\end{proof}

\begin{Lemma}\label{IntersectBorel}
Let $G$ be an almost simple group, whose simple socle $S$ is a group of Lie type, satisfying the conditions of Lemma
{\em \ref{IntersectionOutBorel}}. Let $B=U\leftthreetimes T$ be a Borel subgroup of $\widehat{S}$ and $H= N_G(B)$. Then there exist $x,y,z\in S$, such that $H\cap H^x\cap H^y\cap H^z=1$.
\end{Lemma}

\begin{proof} We use the notations introduced in Lemmas \ref{IntersectionParabolicLeviCenter} and
\ref{IntersectionOutBorel}, in particular $H\leq B\leftthreetimes A$. It is proven in Lemma \ref{IntersectionParabolicLeviCenter}
that there exists $x\in U^{n_0}\leq S$ such that $T^x\cap B=1$. In particular, $B\cap B^{n_0}\cap B^{x^{-1}}=1$. Therefore $H\cap H^{n_0}\cap
H^{x^{-1}}\leq A$ and $A\cap B=1$. By Lemma~\ref{IntersectionOutBorel} there exists $y\in T\cap S=(B\cap B^{n_0})\cap S$ such that $A\cap
A^y=1$. Thus
\begin{equation*}
(H\cap H^{n_0}\cap H^{x^{-1}})\cap  (H\cap H^{n_0}\cap H^{x^{-1}})^y=H\cap H^{n_0}\cap H^{x^{-1}}\cap H^{x^{-1}y}=1,
\end{equation*}
whence the lemma follows.
\end{proof}

\begin{Lemma}\label{ExceptionalOddnotp}
Let $S$ be a simple exceptional group of Lie type over a field of characteristic  $p\not\in\pi$ and $H$ be a solvable
$\pi$-Hall subgroup of $S$. Then one of the followings hold.
\begin{itemize}
 \item[{\em (1)}] There exists a maximal torus $T$ of $S$ such that $H\leq N(S,T)$ and $\vert \pi(N(S,T)/T)\cap
\pi\vert\le1$.
\item[{\em (2)}] $S={}^2G_2(3^{2n+1})$, $\pi\cap\pi(S)=\{2,7\}$, $\vert S\vert_{\{2,7\}}=56$, $H$ is a Frobenius group of order~$56$.
\item[{\em (3)}] $S\in \{G_2(q),F_4(q),E_6^{-\varepsilon}(q), {}^3D_4(q)\}$, where  $\varepsilon\in\{+,-\}$
is chosen so that $q\equiv \varepsilon1\pmod 4${\em;} $2,3\in\pi$, $\pi\cap\pi(S)\subseteq \pi(q-\varepsilon1)$, $H\leq
N(S,T)$, where $T$ is a unique up to conjugation maximal torus such that $N(S,T)$ contains a Sylow
$2$-subgroup of $G$ and $N(S,T)/T$ is a $\{2,3\}$-group. Here $N(S,T):=N_{\overline{G}}(\overline{T})\cap S$, where $T=\overline{T}\cap S$ and $S=O^{p'}(\overline{G}_\sigma)$.
\end{itemize}
\end{Lemma}

\begin{proof} If $2\not\in\pi$ then by \cite[Lemmas 7--14, Theorem~3]{OddHall} statement (1) of the Lemma holds.

If $2\in\pi$ and $3\not\in\pi$, then by \cite[Lemma~5.1 and Theorem~5.2]{Hall3'} (see also
\cite[Theorem~8.9]{VdoRevUspehi}) either statement (1) or statement (2) of the lemma holds.

Finally, if
$2,3\in\pi$, then $S$ is neither a Suzuki group, nor a Ree group (since $p\not\in\pi$). By
\cite[Lemma 7.1--7.6]{VdoRevNumClasses} (see also \cite[Theorem~8.15]{VdoRevUspehi}) we have
$\pi\cap\pi(S)\subseteq \pi(q-\varepsilon1)$, $H\leq
N(S,T)$, where $T$ is a unique up to conjugation maximal torus such that $N(S,T)$ contains a Sylow
$2$-subgroup of $S$ and either $N(S,T)/T$ is a $\{2,3\}$-group, or $N(S,T)/T$ is a Weyl group of the root system of
$S$. Since for root systems  $E_6,E_7,E_8,F_4,G_2$ the Weyl groups are either $\{2,3\}$-groups, or unsolvable, we obtain that if $S\in \{E_6^\varepsilon(q),E_7(q),E_8(q)\}$, then $H$ is unsolvable, whence statement (3) of the lemma.
\end{proof}

\begin{Cor}\label{OrderHnotp}
Let $S$ be a simple exceptional group of Lie type over a field of characteristic $p\not\in\pi$, $S$ is neither a Suzuki group, nor a Ree group, and  $H$ is a solvable $\pi$-Hall subgroup of $S$. Then the following statements hold.
\begin{itemize}
 \item[{\em (1)}] If $S=E_8(q)$, then $\vert H\vert\le (q+1)^8\cdot 2^{14}$.
 \item[{\em (2)}] If $S=E_7(q)$, then $\vert H\vert\le (q+1)^7\cdot 2^{10}$.
 \item[{\em (3)}] If $S=E_6^\epsilon(q)$, then $\vert H\vert\le (q+1)^6\cdot 2^{7}$.
 \item[{\em (4)}] If $S=F_4(q)$, then $\vert H\vert\le (q+1)^4\cdot 2^{7}\cdot 3^2$.
 \item[{\em (5)}] If $S=G_2(q)$, then $\vert H\vert\le (q+1)^2\cdot 12$.
 \item[{\em (6)}] If $S={}^3D_4(q)$, then $\vert H\vert\le \max\{(q^2+q+1)^2,(q+1)^2\cdot 48\}$.
\end{itemize}
\end{Cor}

\section{Proof of the Main Theorem.}

We proceed by considering distinct possible cases for the simple socle $S$ of $G$ and the structure of its
$\pi$-Hall subgroup $H$. If $S$ is either a Suzuki group or a Ree group, then by \cite[Tables 3 and 4]{BurLieSha} it follows that
for every subgroup $H$ of $G$ the inequality $\mathrm{Base}_H(G)\le 3$ holds. So we assume later that $S$ is neither a Suzuki group, nor a Ree group.

\subsection{$S$ is a simple group of Lie type over a field of characteristic $p\in\pi$.}

By Lemma \ref{HallSubgrprop}, $H\cap \widehat{S}\in\mathrm{Hall}_\pi(\widehat{S})$, so in this case for $H\cap \widehat{S}$ Lemma \ref{SolvHallcharacteristic} holds. Assume first that
$H\cap \widehat{S}$ lies in a Borel subgroup of $\widehat{S}$. If $S$ is a nontwisted group of Lie type over a field of order two, then $H$ is a  $2$-group. By \cite{Zen1} the inequality $\mathrm{Base}_H(G)\le 3$ holds. If
$S=D_4(3)$, then $H$ is a $3$-group. By \cite{Zen1} the inequality  $\mathrm{Base}_H(G)\le 3$ holds. Assume that  $S$ is not a nontwisted group of Lie type over a field of two elements, and $S\not\simeq D_4(3)$. Then $H\leq
N_G(U)=N_G(B)$ and by Lemma \ref{IntersectBorel} the inequality $\mathrm{Base}_H(G)\le 4$ holds. If one of statements (1)--(3) of Lemma \ref{SolvHallcharacteristic} is satisfied, then $S$ is a classical group and calculations by using
\cite{GAP} show that in any case $\mathrm{Base}_H(G)\le 5$ and~${\mathrm{Reg}_H(G,5)\ge 5}$.

\subsection{$S$ is a simple exceptional group of Lie type over a field of characteristic $p\not\in\pi$.}

Asuume that $S=E_8(q)$. We use Lemma \ref{BoundfprBase}. If $x$ is a unipotent element, then $x^G\cap
H=\varnothing$. If  $x$ is a semisimple element from $G=\widehat{G}$, then by \cite[Table~2]{Der2} it follows that the maximum
of orders of centralizers of semisimple elements in $E_8(q)$ is not greater than
\begin{equation*}
q^{64}(q^{18}-1)(q^{14}-1)(q^{12}-1)(q^{10}-1)(q^8-1)(q^6-1)(q^2-1)^2,
\end{equation*}
whence $\vert x^G\vert>q^{112}$. Clealy, the inequality $\vert x^G\vert>q^{112}$ holds in case, when $x$ is a field automorphism. So for $c=2$
we obtain
\begin{equation*}
\widehat{Q}(G,2)\leq ((q+1)^8\cdot 2^{14})^2/(q^{112})<1
\end{equation*}
for every $q\ge 2$. Hence,
$\mathrm{Base}_H(G)\le 2$.

Assume that $G=E_7(q)$. We again use Lemma \ref{BoundfprBase}. If $x$ is a unipotent element, then $x^G\cap
H=\varnothing$. If $x$ is a semisimple element from  $\widehat{G}$, then by \cite[Table~1]{Der2} it follows that the maximum of orders of centralizers
of semisimple elements in $E_7(q)$ is not greater than
\begin{equation*}
 q^{31}(q^2-1)^2(q^4-1)(q^6-1)^2(q^8-1)(q^{10}-1),
\end{equation*}
whence $\vert x^G\vert>(1/2)q^{64}$. Clearly, the inequality
$\vert x^G\vert>(1/2)q^{64}$ hols in case, when $x$ is a field automorphism. So for $c=2$
we obtain
\begin{equation*}
 \widehat{Q}(G,2)\leq ((q+1)^7\cdot 2^{20})^2\cdot 2/(q^{64})<1
\end{equation*}
for every $q\ge 2$. Hence $\mathrm{Base}_H(G)\le 2$.

Assume that $G=E_6^\epsilon(q)$. As above, we obtain that $x$ is either a semisimple element from  $\widehat{G}$, or does not lie in $\widehat{G}$. If $x$ is a semisimple element, then by \cite[Table 1 and Case $E_6(q)$]{Der1} it follows that the maximum of orders of
of centralizers of semisimple elements in  $E_6^\epsilon(q)$ is not greater than
\begin{equation*}
q^{20}(q-\epsilon1)(q^2-1)(q^4-1)(q^6-1)(q^8-1)(q^5-\epsilon1),
\end{equation*}
whence $\vert x^G\vert>\frac{1}{3} q^{30}$. Clearly, the inequality $\vert x^G\vert>\frac{1}{3} q^{30}$ holds in case, when   $x$ is either a field, or a graph-field automorphism. If $x$ is a graph automorphism, then
\begin{equation*}
\vert x^G\vert =\vert E_6^\epsilon\vert/\vert F_4(q)\vert\ge \frac{1}{3}q^{12}(q^5-1)(q^9-1).
\end{equation*}
So for $c=4$ we obtain
\begin{equation*}
\widehat{Q}(G,2)\leq \frac{(q+1)^{24}\cdot 2^{28}\cdot 3^3}{q^{36}\cdot
(q^5-1)^3\cdot(q^9-1)^3}<1
\end{equation*}
for every $q\ge 2$. Hence
$\mathrm{Base}_H(G)\le 4$.

Assume that $G=F_4(q)$. Again we may assume that $x$ either is a semisimple element from $G=\widehat{G}$ or does not lie in
$\widehat{G}$. If $x$ is a semisimple element, then by \cite[Table 2]{Der1} it follows that the maximum of orders of centralizers
of semisimple elements in  $F_4(q)$ is not greater than
\begin{equation*}
q^{16}(q^2-1)(q^4-1)(q^6-1)(q^8-1),
\end{equation*}
whence $\vert x^G\vert > q^{16}$. Clearly, the inequality  $\vert x^G\vert > q^{16}$ holds for every $x$ not lying in
$G$. So for $c=4$ we obtain
\begin{equation*}
\widehat{Q}(G,2)\leq \frac{(q+1)^{16}\cdot 2^{28}\cdot 3^8}{q^{48}}<1
\end{equation*}
for every $q\ge 3$. So for $q\ge 3$ the inequality $\mathrm{Base}_H(G)\le 4$ holds. If $q=2$, then in view of the condition
$p\not\in\pi$ we obtain that the order $\vert H\vert$ is odd. By \cite[Lemma~8]{OddHall} we obtain that either  $H$ is a Sylow $3$-subgroup of $G$, or $H$ is abelian.
Hence the inequality  $\mathrm{Base}_H(G)\le 3$ holds: in the first case by \cite{Zen1}, and in the second case by Lemma~\ref{ZenkovAbelian}.

Assume that $G=G_2(q)$. As above we may assume that $x$  either is a semisimple element from $G=\widehat{G}$, or does
not lie in $\widehat{G}$. If $x$ is semisimple, then by \cite[Table 4]{Der1} it follows that the maximum of orders of centralizers of semisimple elements in $F_4(q)$ is not greater than
\begin{equation*}
q^{2}(q^2-1)(q^3+1),
\end{equation*}
whence $\vert x^G\vert\ge q^4(q^3-1)$. Clearly, the inequality $\vert x^G\vert\ge q^4(q^3-1)$ holds for every $x$ not lying in $G$. So for $c=4$ we obtain
\begin{equation*}
\widehat{Q}(G,2)\leq \frac{(q+1)^8\cdot 12^4}{q^{12}\cdot(q^3-1)^3}<1
\end{equation*}
for every $q\ge 3$. Hence for $q\ge 3$ the inequality $\mathrm{Base}_H(G)\le 4$ holds. If $q=2$, then by the condition
$p\not\in\pi$ we obtain that the order $\vert H\vert$ is odd. By \cite[Lemma~7]{OddHall} we obtain that either $H$ is a Sylow $3$-subgroup of $G$, or $H$ is abelian. So the inequality $\mathrm{Base}_H(G)\le 3$ holds: in the first case by \cite{Zen1}, and in the second case by Lemma~\ref{ZenkovAbelian}.

If $G={}^3D_4(q)$, then by \cite[Table~7]{Der1} it is easy to get the bound  $\vert x^G\vert>q^{16}$. Using the bound we obtain that for
$q\ge 2$ the inequality  $\mathrm{Base}_H(G)\le 4$ holds. The Main Theorem is proven.

Notice that for the case $p\in\pi$ we also prove the following

\begin{Theo}\label{Liep}
Let $G$ be a finite almost simple group, whose simple socle is isomorphic to a group of Lie type over  a field of characteristic
$p\in\pi$. Assume that $H$ is a solvable $\pi$-Hall subgroup of $G$. Then the inequalities
$\mathrm{Base}_H(G)\leq 5$ and $\mathrm{Reg}_H(G,5)\geq 5$ hold.
\end{Theo}


\end{document}